\documentclass[11pt]{article}
\usepackage{amsfonts}
\usepackage{amsfonts}
\usepackage{amsfonts}
\usepackage{mathrsfs}
\usepackage{bm}
\usepackage{amssymb,amsmath} 
 \allowdisplaybreaks

\catcode`!=11
\let\!int\int \def\int{\displaystyle\!int}
\let\!lim\lim \def\lim{\displaystyle\!lim}
\let\!sum\sum \def\sum{\displaystyle\!sum}
\let\!sup\sup \def\sup{\displaystyle\!sup}
\let\!inf\inf \def\inf{\displaystyle\!inf}
\let\!cap\cap \def\cap{\displaystyle\!cap}
\let\!max\max \def\max{\displaystyle\!max}
\let\!min\min \def\min{\displaystyle\!min}
\let\!frac\frac \def\frac{\displaystyle\!frac}
\catcode`!=12

\let\oldsection\section
\renewcommand\section{\setcounter{equation}{0}\oldsection}

\allowdisplaybreaks
\def\pf{\it{Proof.}\rm\quad}

\newtheorem{thm}{Theorem}[section]
\newtheorem{lem}[thm]{Lemma}
\newtheorem{cor}[thm]{Corollary}

\begin{document}
\title {\bf Some Evaluation of Quadratic Euler Sums}
\author{
{Ce Xu\thanks{Corresponding author. Email: 15959259051@163.com}}\\[1mm]
\small School of Mathematical Sciences, Xiamen University\\
\small Xiamen
361005, P.R. China}

\date{}
\maketitle \noindent{\bf Abstract }  In this paper, we obtain some formulas for double nonlinear Euler sums involving harmonic numbers and alternating harmonic numbers. By using these formulas, we give new closed form
sums of several quadratic Euler series through Riemann zeta values, polylogarithm functions and linear sums. Furthermore, some relationships between Euler sums and integrals of polylogarithm functions are established.
\\[2mm]
\\[2mm]

\section{Introduction}
Throughout this article we will use the following definitions and notations. In this paper, the generalized harmonic numbers and alternating harmonic numbers is defined by
$$H_n=\sum\limits_{j=1}^n\frac {1}{j},\ \zeta_n(k)=\sum\limits_{j=1}^n\frac {1}{j^k} ,\ L_{n}(k)=\sum\limits_{j=1}^n\frac{(-1)^{j-1}}{j^k},\ 1\leq k \in Z.\eqno(1.1)$$
For a pair $(p,q)$ of positive integers with $q \geq 2$, the classical double linear Euler sum is defined by
$$S_{p,q}=\sum\limits_{n = 1}^\infty  {\frac{1}{{{n^q}}}} \sum\limits_{k = 1}^n {\frac{1}{{{k^p}}}}, \eqno(1.2)$$
The number $w=p+q$ is defined as the weight of $S_{p,q}$. The evaluation of $S_{p,q}$ in terms of values of Riemann zeta function at positive integers is known when $p=1,\ p=q,\ (p,q)=(2,4),(4,2)$ or $p+q$ is odd.
In 1742, Goldbach proposed to Euler the problem of expressing the $S_{p,q}$ in terms of values at positive integers of the Riemann zeta function $\zeta(s)$. The Riemann zeta function and  alternating Riemann zeta function are defined respectively by $$\zeta(s)=\sum\limits_{n = 1}^\infty {\frac {1}{n^{s}}},\Re(s)>1, $$
and
$$\bar \zeta \left( s \right) = \sum\limits_{n = 1}^\infty  {\frac{{{{\left( { - 1} \right)}^{n - 1}}}}{{{n^s}}}} ,\;{\mathop{\Re}\nolimits} \left( s \right) \ge 1.$$
Euler proved that the double linear sums are reducible to zeta values whenever $a + b$ is less
than 7 or when $a + b$ is odd and less than 13. He conjectured that the double linear
sums would be reducible to zeta values when $p + q$ is odd, and even gave what he hoped to obtain the general
formula. In [3], D. Borwein, J.M. Borwein and R. Girgensohn proved conjecture and formula,
and in [2], D.H. Bailey, J.M. Borwein and R. Girgensohn conjectured that the double linear sums when $p + q > 7,p + q$ is even, are not reducible.
In [12], Philippe Flajolet and Bruno Salvy gave a general formula for odd weight $p+q$,
\begin{align*}
\sum\limits_{n = 1}^\infty  {\frac{{{\zeta _n}\left( p \right)}}{{{n^q}}}}
 &= \zeta \left( m \right)\left( {\frac{1}{2} - \frac{{{{\left( { - 1} \right)}^p}}}{2}\left( {\begin{array}{*{20}{c}}
   {m - 1}  \\
   p  \\
\end{array}} \right) - \frac{{{{\left( { - 1} \right)}^p}}}{2}\left( {\begin{array}{*{20}{c}}
   {m - 1}  \\
   q  \\
\end{array}} \right)} \right)
\nonumber \\
           &\quad  + {\left( { - 1} \right)^p}\sum\limits_{k = 1}^{\left[ {p/2} \right]} {\left( {\begin{array}{*{20}{c}}
   {m - 2k - 1}  \\
   {q - 1}  \\
\end{array}} \right)\zeta \left( {2k} \right)\zeta \left( {m - 2k} \right)} + \frac{{1 - {{\left( { - 1} \right)}^p}}}{2}\zeta \left( p \right)\zeta \left( q \right)
 \nonumber \\
           &\quad + {\left( { - 1} \right)^p}\sum\limits_{k = 1}^{\left[ {q/2} \right]} {\left( {\begin{array}{*{20}{c}}
   {m - 2k - 1}  \\
   {p - 1}  \\
\end{array}} \right)\zeta \left( {2k} \right)\zeta \left( {m - 2k} \right)},
\end{align*}
where $\zeta(1)$ should be interpreted as 0 wherever it occurs.\\
Let $\pi  = \left( {{\pi _1}, \ldots ,{\pi _k}} \right)$ be a partition of integer $p$ and $p = {\pi _1} +  \cdots  + {\pi _k}$ with ${\pi _1} \le {\pi _2} \le  \cdots  \le {\pi _k}$. The classical double nonlinear Euler sum of index $\pi,q$ is defined as follows (see [12])
\[{S_{\pi ,q}} = \sum\limits_{n = 1}^\infty  {\frac{{{\zeta _n}\left( {{\pi _1}} \right){\zeta _n}\left( {{\pi _2}} \right) \cdots {\zeta _n}\left( {{\pi _k}} \right)}}{{{n^q}}}} ,\tag{1.3}\]
where the quantity ${\pi _1} +  \cdots  + {\pi _k} + q$ is called the weight, the quantity $k$ called the degree. As usual, repeated summands in partitions are indicated by powers, so that for instance
\[{S_{{1^2}{2^3}4,q}} = {S_{112224,q}} = \sum\limits_{n = 1}^\infty  {\frac{{H_n^2\zeta _n^3\left( 2 \right){\zeta _n}\left( 4 \right)}}{{{n^q}}}}. \]
The general Euler sums are defined by the series (see [19-21])
$$\sum\limits_{n = 1}^\infty  {\frac{{\prod\limits_{i = 1}^{{m_1}} {\zeta _n^{{q_{_i}}}\left( {{k_i}} \right)\prod\limits_{j = 1}^{{m_2}} {L_n^{{l_j}}\left( {{h_j}} \right)} } }}{{{n^p}}}} ,\ \sum\limits_{n = 1}^\infty  {\frac{{\prod\limits_{i = 1}^{{m_1}} {\zeta _n^{{q_{_i}}}\left( {{k_i}} \right)\prod\limits_{j = 1}^{{m_2}} {L_n^{{l_j}}\left( {{h_j}} \right)} {{\left( { - 1} \right)}^{n - 1}}} }}{{{n^p}}}}, \eqno(1.4)$$
where $p (p>1),m_{1},m_{2},q_{i},k_{i},h_{j},l_{j}$ are positive integer. If $\sum\limits_{i=1}^{m_{1}}(k_{i}q_{i})+\sum\limits_{j=1}^{m_{2}}(h_{j}l_{j})+p=C$($C$ is a positive integer), then we call the identity $C$th-order Euler sums.\\
In [13], Philippe Flajolet and Bruno Salvy gave explicit reductions to zeta values and logarithm for all linear sums
\[\sum\limits_{n = 1}^\infty  {\frac{{{\zeta _n}\left( p \right)}}{{{n^q}}}} ,\;\sum\limits_{n = 1}^\infty  {\frac{{{L_n}\left( p \right)}}{{{n^q}}}} ,\;\sum\limits_{n = 1}^\infty  {\frac{{{\zeta _n}\left( p \right)}}{{{n^q}}}} {\left( { - 1} \right)^{n - 1}},\;\sum\limits_{n = 1}^\infty  {\frac{{{L_n}\left( p \right)}}{{{n^q}}}{{\left( { - 1} \right)}^{n - 1}}} \]
when $p+q$ is an odd weight. The evaluation of linear sums in
terms of values of the Riemann zeta function and polylogarithm function at positive integers is known when
 $(p,q) = (1,3), (2,2)$, or $p + q$ is odd. For example
\begin{align*}
&\sum\limits_{n = 1}^\infty  {\frac{{{H_n}{}}}{{{n^3}}}}{\left( { - 1} \right)}^{n - 1}  =  - 2{ Li}{_4}\left( {\frac{1}{2}} \right) + \frac{{11{}}}{{4}} \zeta(4) + \frac{{{1}}}{{2}}\zeta(2){\ln ^2}2 - \frac{1}{{12}}{\ln ^4}2 - \frac{7}{4}\zeta \left( 3 \right)\ln 2,\\
&\sum\limits_{n = 1}^\infty  {\frac{{{L_n}\left( 1 \right)}}{{{n^3}}}{{\left( { - 1} \right)}^{n - 1}} = \frac{3}{2}\zeta \left( 4 \right)}  + \frac{1}{2}\zeta \left( 2 \right){\ln ^2}2 - \frac{1}{{12}}{\ln ^4}2 - 2{ Li}{_4}\left( {\frac{1}{2}} \right),\\
&\sum\limits_{n = 1}^\infty  {\frac{{{L_n}\left( 2 \right)}}{{{n^2}}} = \frac{{85}}{{16}}\zeta \left( 4 \right)}  - 4{ Li}{_4}\left( {\frac{1}{2}} \right) + \zeta \left( 2 \right){\ln ^2}2 - \frac{1}{6}{\ln ^4}2 - \frac{7}{2}\zeta \left( 3 \right)\ln 2,\\
&\sum\limits_{n = 1}^\infty  {\frac{{{}{\zeta _n}\left( 2 \right)}}{{{n^2}}}} {\left( { - 1} \right)}^{n - 1}=-\frac{51}{16}\zeta(4) +4 {Li}_4 \left(\frac{1}{2} \right)+\frac{7}{2}\ln2 \zeta(3)-\zeta(2) \ln^2 2+\frac{\ln^4 2}{6}.
\end{align*}
In many other cases we are not able to obtain a formula for the Euler sum constant
explicitly in terms of values of the Riemann zeta, logarithm and polylogarithm functions,
but we are able to obtain relations involving two of more Euler sum constants of the same
degree.
In [21], we proved that the quadratic double sums
\[\sum\limits_{n = 1}^\infty  {\frac{{H_n^2}}{{{n^{2m}}}}} {\left( { - 1} \right)^{n - 1}},\sum\limits_{n = 1}^\infty  {\frac{{L_n^2\left( 1 \right)}}{{{n^{2m}}}}} {\left( { - 1} \right)^{n - 1}},\sum\limits_{n = 1}^\infty  {\frac{{{H_n}{L_n}\left( 1 \right)}}{{{n^{2m}}}}} {\left( { - 1} \right)^{n - 1}},\sum\limits_{n = 1}^\infty  {\frac{{{H_n}{L_n}\left( 1 \right)}}{{{n^{2m}}}}} ,\sum\limits_{n = 1}^\infty  {\frac{{L_n^2\left( 1 \right)}}{{{n^{2m}}}}} \]
are reducible to polynomials zeta values and to linear sums, $m$ is a positive integer.
The relationship between the
values of the Riemann zeta function and Euler sums has been studied by many authors, for example see [2-4,6-21].\\
The main purpose of this paper is to evaluate some quadratic Euler sums which involving harmonic numbers and alternating harmonic numbers, either linearly or nonlinearly. In this paper, we will prove that all quadratic double sums
\[\sum\limits_{n = 1}^\infty  {\frac{{{\zeta _n}\left( p \right){\zeta _n}\left( {p + 2m + 1} \right)}}{n}} {\left( { - 1} \right)^{n - 1}},\sum\limits_{n = 1}^\infty  {\frac{{{L_n}\left( p \right){L_n}\left( {p + 2m + 1} \right)}}{n}} {\left( { - 1} \right)^{n - 1}},\ 2 \le p \in Z,\;0 \le m \in Z\]
are reducible to linear sums and to polynomials zeta values. In the same way, we also obtain that, for $2 \le p \in Z,\;0 \le m \in Z$, the following expression
$$\sum\limits_{n = 1}^\infty  {\left\{ {\frac{{{H_n}{\zeta _n}\left( {p + 2m + 1} \right)}}{{{n^p}}} + \frac{{{H_n}{\zeta _n}\left( p \right)}}{{{n^{p + 2m + 1}}}}} \right\}} $$
and
$$\sum\limits_{n = 1}^\infty  {\left\{ {\frac{{{H_n}{\zeta _n}\left( {p + 2m + 2} \right)}}{{{n^p}}} - \frac{{{H_n}{\zeta _n}\left( p \right)}}{{{n^{p + 2m + 2}}}}} \right\}} $$ are reducible to linear sums.

\section{Main Theorems and Proof}
In this section, we will establish some explicit relationships which
involve Euler sums and integrals of polylogarithm functions. The polylogarithm function defined as follows
\[L{i_p}\left( x \right) = \sum\limits_{n = 1}^\infty  {\frac{{{x^n}}}{{{n^p}}}}, \Re(p)>1,\ \left| x \right| \le 1 .\]
when $x$ takes $1$ and $ -1$, then the function ${L{i_p}\left( {x} \right)}$ are reducible to Riemann zeta function and alternating Riemann zeta function, respectively.
\begin{lem}\label{lem 2.1}
Let $s,t\ge 1$ be integers. Then the product $L{i_s}\left( x \right)L{i_t}\left( x \right)$ are reducible to Euler sums function:
\begin{align*}
 L{i_s}\left( x \right)L{i_t}\left( x \right)  &= \sum\limits_{j = 1}^s {A_j^{\left( {s,t} \right)}} \sum\limits_{n = 1}^\infty  {\frac{{{\zeta _n}\left( j \right)}}{{{n^{s + t - j}}}}} {x^n} + \sum\limits_{j = 1}^t {B_j^{\left( {s,t} \right)}} \sum\limits_{n = 1}^\infty  {\frac{{{\zeta _n}\left( j \right)}}{{{n^{s + t - j}}}}} {x^n}
\nonumber \\ &\quad - \left( {\sum\limits_{j = 1}^s {A_j^{\left( {s,t} \right)}}  + \sum\limits_{j = 1}^t {B_j^{\left( {s,t} \right)}} } \right)L{i_{s + t}}\left( x \right).\tag{2.1}
\end{align*}
where $A_j^{\left( {s,t} \right)} = \left( {\begin{array}{*{20}{c}}
   {s + t - j - 1}  \\
   {s - j}  \\
\end{array}} \right),B_j^{\left( {s,t} \right)} = \left( {\begin{array}{*{20}{c}}
   {s + t - j - 1}  \\
   {t - j}  \\
\end{array}} \right).$
\end{lem}
\pf
The left hand of equation (2.1) equals
\[L{i_s}\left( x \right)L{i_t}\left( x \right) = \sum\limits_{n = 1}^\infty  {\sum\limits_{k = 1}^n {\frac{{{x^{n + 1}}}}{{{k^s}{{\left( {n - k + 1} \right)}^t}}}} } .\tag{2.2}\]
By the formula
\[\frac{1}{{{x^s}{{\left( {1 - x} \right)}^t}}} = \sum\limits_{j = 1}^s {\frac{{A_j^{\left( {s,t} \right)}}}{{{x^j}}}}  + \sum\limits_{j = 1}^t {\frac{{B_j^{\left( {s,t} \right)}}}{{{{\left( {1 - x} \right)}^j}}}} ,\;\;s,t \ge 0,s + t \ge 1.\tag{2.3}\]
 we can obtain $(2.1)$.
\begin{thm}\label{lem 2.2}\ \ {\sl Let $p,q \ge 1$ be integers and $x \in \left[ { - 1,1} \right)$. Then the following identity holds}
\begin{align*}
 \int\limits_0^x {\frac{{L{i_p}\left( t \right)L{i_q}\left( t \right)}}{t}dt} &= \sum\limits_{i = 1}^{q - 1} {{{\left( { - 1} \right)}^{i - 1}}L{i_{p + i}}\left( x \right)L{i_{q + 1 - i}}\left( x \right)}  + {\left( { - 1} \right)^q}\ln \left( {1 - x} \right)\left( {L{i_{p + q}}\left( x \right) - \zeta \left( {p + q} \right)} \right)
\nonumber \\ &\quad  - {\left( { - 1} \right)^q}\sum\limits_{n = 1}^\infty  {\frac{1}{{{n^{p + q}}}}\left( {\sum\limits_{k = 1}^n {\frac{{{x^k}}}{k}} } \right)} . \tag{2.4}
\end{align*}
\end{thm}
\pf Let \[{I_{p,q}}\left( x \right) = \int\limits_0^x {\frac{{L{i_p}\left( t \right)L{i_q}\left( t \right)}}{t}dt} .\tag{2.5}\]
By the definition of polylogarithm function, we can verify that
\[\int\limits_0^x {\frac{{L{i_p}\left( t \right)L{i_q}\left( t \right)}}{t}dt}  = \sum\limits_{n = 1}^\infty  {\frac{1}{{{n^p}}}} \int\limits_0^x {{t^{n - 1}}L{i_q}\left( t \right)dt} .\tag{2.6}\]
is hold.
Using integration by parts we have
\[\int\limits_0^x {{t^{n - 1}}L{i_q}\left( t \right)dt}  = \sum\limits_{i = 1}^{q - 1} {{{\left( { - 1} \right)}^{i - 1}}\frac{{{x^n}}}{{{n^i}}}L{i_{q + 1 - i}}\left( x \right)}  + \frac{{{{\left( { - 1} \right)}^q}}}{{{n^q}}}\ln \left( {1 - x} \right)\left( {{x^n} - 1} \right) - \frac{{{{\left( { - 1} \right)}^q}}}{{{n^q}}}\left( {\sum\limits_{k = 1}^n {\frac{{{x^k}}}{k}} } \right).\tag{2.7}\]
Combining (2.6) and (2.7), we obtain (2.4) holds.\\
Euler gave the following formula in 1775 (see [2, 12])
\[\sum\limits_{n = 1}^\infty  {\frac{{{H_n}}}{{{n^k}}}}  = \frac{1}{2}\left\{ {\left( {k + 2} \right)\zeta \left( {k + 1} \right) - \sum\limits_{i = 1}^{k - 2} {\zeta \left( {k - i} \right)\zeta \left( {i + 1} \right)} } \right\}.\tag{2.8}\]
Letting $x\rightarrow1$ in Theorem 2.1, we can get the following results
\begin{cor}(\cite{F2005})
Let $p,q$ be integers with $p,q \ge 1$. Then the integrals ${I_{p,q}}\left( 1 \right)$ reduce to polynomials in zeta values:
\begin{align*}
 \int\limits_0^1 {\frac{{L{i_p}\left( x \right)L{i_q}\left( x \right)}}{x}} dx &= \sum\limits_{i = 1}^{q - 1} {{{\left( { - 1} \right)}^{i - 1}}} \zeta \left( {q + 1 - i} \right)\zeta \left( {p + i} \right) + {\left( { - 1} \right)^{q - 1}}\left( {1 + \frac{{p + q}}{2}} \right)\zeta \left( {p + q + 1} \right)
\nonumber \\ &\quad - \frac{1}{2}\sum\limits_{k = 1}^{p + q - 2} {\zeta \left( {k + 1} \right)} \zeta \left( {p+q - k} \right).
\end{align*}
\end{cor}
\begin{thm} Let $s,t$ be integers with $s,t\ge 1$, we have
\begin{align*}
 \int\limits_0^1 {\frac{{Li_s^2\left( x \right)L{i_t}\left( x \right)}}{x}} dx
 &=2\sum\limits_{j = 1}^s {A_j^{\left( {s,s} \right)}\left\{ {\sum\limits_{i = 1}^{t - 1} {{{\left( { - 1} \right)}^{i - 1}}\zeta \left( {t + 1 - i} \right)\sum\limits_{n = 1}^\infty  {\frac{{{\zeta _n}\left( j \right)}}{{{n^{2s + i - j}}}}} }  + {{\left( { - 1} \right)}^{t - 1}}\sum\limits_{n = 1}^\infty  {\frac{{{H_n}{\zeta _n}\left( j \right)}}{{{n^{2s + t - j}}}}} } \right\}}
  \nonumber \\
           &\quad - 2\sum\limits_{j = 1}^s {A_j^{\left( {s,s} \right)}} I\left( {2s,t} \right),\tag{2.9}
\end{align*}
\begin{align*}
 \int\limits_0^1 {\frac{{Li_s^2\left( x \right)L{i_t}\left( x \right)}}{x}} dx
 &= {\left( { - 1} \right)^{s - 1}}\sum\limits_{n = 1}^\infty  {\left\{ {\sum\limits_{j = 1}^s {A_j^{\left( {s,t} \right)}} \frac{{{H_n}{\zeta _n}\left( j \right)}}{{{n^{2s + t - j}}}} + \sum\limits_{j = 1}^t {B_j^{\left( {s,t} \right)}} \frac{{{H_n}{\zeta _n}\left( j \right)}}{{{n^{2s + t - j}}}}} \right\}}
\nonumber \\
& \quad +\sum\limits_{i = 1}^{s - 1} {{{\left( { - 1} \right)}^{i - 1}}\zeta \left( {s + 1 - i} \right)\sum\limits_{n = 1}^\infty  {\left\{ {\sum\limits_{j = 1}^s {A_j^{\left( {s,t} \right)}} \frac{{{\zeta _n}\left( j \right)}}{{{n^{s + t + i - j}}}} + \sum\limits_{j = 1}^t {B_j^{\left( {s,t} \right)}} \frac{{{\zeta _n}\left( j \right)}}{{{n^{s + t + i - j}}}}} \right\}} }
\nonumber \\
& \quad - \left\{ {\sum\limits_{j = 1}^s {A_j^{\left( {s,t} \right)}}  + \sum\limits_{j = 1}^t {B_j^{\left( {s,t} \right)}} } \right\}I\left( {s,s + t} \right)
.\tag{2.10}
\end{align*}
\end{thm}
\pf Let $s=t$ in $(2.1)$, we have
\begin{align*}
Li_s^2\left( x \right) = 2\sum\limits_{j = 1}^s {A_j^{\left( {s,s} \right)}} \sum\limits_{n = 1}^\infty  {\frac{{{\zeta _n}\left( j \right)}}{{{n^{2s - j}}}}} {x^n} - 2\sum\limits_{j = 1}^s {A_j^{\left( {s,s} \right)}} L{i_{2s}}\left( x \right).\tag{2.11}
\end{align*}
Multiplying $\frac{{L{i_t}\left( x \right)}}{x}$ to the equation (2.11) and integrating over $(0,1)$, by virtue of $(2.7)$, we obtain $(2.9)$.\\
Similarly, multiplying $\frac{{L{i_s}\left( x \right)}}{x}$ in both sides of $(2.1)$, and integrating from $x=0$ to $x=1$ with $(2.7)$. By a simple calculation, we obtain formula $(2.10)$.\\
\begin{thm}
For $p\geq 2,m\geq 0$ and $p,m\in Z$,$x \in \left[ { - 1,1} \right)$
, we have
\begin{align*}
&{\left( { - 1} \right)^p}\sum\limits_{n = 1}^\infty  {\left\{ {\frac{{{\zeta _n}\left( {p + 2m + 1} \right)}}{{{n^p}}} + \frac{{{\zeta _n}\left( p \right)}}{{{n^{p + 2m + 1}}}}} \right\}\left( {\sum\limits_{k = 1}^n {\frac{{{x^k}}}{k}} } \right)}  \\
& =\left( {p + 2m + 2} \right){I_{p - 1,p + 2m + 2}}\left( x \right) + \left( {2m + 1} \right){I_{p,p + 2m + 1}}\left( x \right) - \left( {p + 1} \right){I_{p + 1,p + 2m}}\left( x \right) \\
&\quad  + {\left( { - 1} \right)^p}\ln \left( {1 - x} \right)\sum\limits_{n = 1}^\infty  {\left\{ {\frac{{{\zeta _n}\left( {p + 2m + 1} \right)}}{{{n^p}}} + \frac{{{\zeta _n}\left( p \right)}}{{{n^{p + 2m + 1}}}}} \right\}\left( {{x^n} - 1} \right)}  \\
&\quad +\sum\limits_{i = 1}^{p + 2m} {{{\left( { - 1} \right)}^{i - 1}}L{i_{p + 2m + 2 - i}}\left( x \right)\sum\limits_{n = 1}^\infty  {\left\{ {\sum\limits_{j = 2}^{p - 1} {\frac{{{\zeta _n}\left( j \right)}}{{{n^{p + i - j}}}}}  + 2\frac{{{H_n}}}{{{n^{p + i - 1}}}}} \right\}{x^n}} }  \\
&\quad  + \sum\limits_{i = 1}^{p + 2m - 1} {{{\left( { - 1} \right)}^{i - 1}}L{i_{p + 2m + 1 - i}}\left( x \right)\sum\limits_{n = 1}^\infty  {\left\{ {\sum\limits_{j = 2}^p {\frac{{{\zeta _n}\left( j \right)}}{{{n^{p + 1 + i - j}}}}}  + 2\frac{{{H_n}}}{{{n^{p + i}}}}} \right\}{x^n}} } \\
&\quad - \sum\limits_{i = 1}^{p - 2} {{{\left( { - 1} \right)}^{i - 1}}L{i_{p - i}}\left( x \right)\sum\limits_{n = 1}^\infty  {\left\{ {\sum\limits_{j = 2}^{p + 2m + 1} {\frac{{{\zeta _n}\left( j \right)}}{{{n^{p + 2m + 2 + i - j}}}}}  + 2\frac{{{H_n}}}{{{n^{p + 2m + 1 + i}}}}} \right\}{x^n}} } \\
&\quad  - \sum\limits_{i = 1}^{p - 1} {{{\left( { - 1} \right)}^{i - 1}}L{i_{p + 1 - i}}\left( x \right)\sum\limits_{n = 1}^\infty  {\left\{ {\sum\limits_{j = 2}^{p + 2m} {\frac{{{\zeta _n}\left( j \right)}}{{{n^{p + 2m + 1 + i - j}}}}}  + 2\frac{{{H_n}}}{{{n^{p + 2m + i}}}}} \right\}{x^n}} }  .\tag{2.12}
\end{align*}
where the integral ${I_{p,q}}\left( x \right)$ is defined by (2.5).
\end{thm}
\pf We first consider the following integral
$\int\limits_0^x {\frac{{\ln \left( {1 - t} \right)L{i_p}\left( t \right)L{i_{p + 2m}}\left( t \right)}}{t}} dt,\;x \in \left[ { - 1,1} \right).$ In (2.1), taking $s=1,t=p$, we have that
\[\ln \left( {1 - t} \right)L{i_p}\left( t \right) = \left( {p + 1} \right)L{i_{p + 1}}\left( t \right) - \sum\limits_{n = 1}^\infty  {\left\{ {\sum\limits_{j = 2}^p {\frac{{{\zeta _n}\left( j \right)}}{{{n^{p + 1 - j}}}}}  + 2\frac{{{H_n}}}{{{n^p}}}} \right\}{t^n}} ,x \in \left[ { - 1,1} \right).\tag{2.13}\]
By virtue of (2.7) and (2.13), we obtain
\begin{align*}
&\int\limits_0^x {\frac{{\ln \left( {1 - t} \right)L{i_p}\left( t \right)L{i_{p + 2m}}\left( t \right)}}{t}} dt \\
& =\left( {p + 1} \right){I_{p + 1,p + 2m}}\left( x \right) - \sum\limits_{i = 1}^{p + 2m - 1} {{{\left( { - 1} \right)}^{i - 1}}L{i_{p + 2m + 1 - i}}\left( x \right)\sum\limits_{n = 1}^\infty  {\left\{ {\sum\limits_{j = 2}^p {\frac{{{\zeta _n}\left( j \right)}}{{{n^{p + 1 + i - j}}}}}  + 2\frac{{{H_n}}}{{{n^{p + i}}}}} \right\}{x^n}} }  \\
&\quad   + {\left( { - 1} \right)^{p - 1}}\ln \left( {1 - x} \right)\sum\limits_{n = 1}^\infty  {\left\{ {\sum\limits_{j = 2}^p {\frac{{{\zeta _n}\left( j \right)}}{{{n^{2p + 2m + 1 - j}}}}}  + 2\frac{{{H_n}}}{{{n^{2p + 2m}}}}} \right\}\left( {{x^n} - 1} \right)}  \\
&\quad  - {\left( { - 1} \right)^{p - 1}}\sum\limits_{n = 1}^\infty  {\left\{ {\sum\limits_{j = 2}^p {\frac{{{\zeta _n}\left( j \right)}}{{{n^{2p + 2m + 1 - j}}}}}  + 2\frac{{{H_n}}}{{{n^{2p + 2m}}}}} \right\}\left( {\sum\limits_{k = 1}^n {\frac{{{x^k}}}{k}} } \right)}   \\
&= \left( {p + 2m + 1} \right){I_{p,p + 2m + 1}}\left( x \right) - \sum\limits_{i = 1}^{p - 1} {{{\left( { - 1} \right)}^{i - 1}}L{i_{p + 1 - i}}\left( x \right)\sum\limits_{n = 1}^\infty  {\left\{ {\sum\limits_{j = 2}^{p + 2m} {\frac{{{\zeta _n}\left( j \right)}}{{{n^{p + 2m + 1 + i - j}}}}}  + 2\frac{{{H_n}}}{{{n^{p + 2m + i}}}}} \right\}{x^n}} }  \\
&\quad  + {\left( { - 1} \right)^{p - 1}}\ln \left( {1 - x} \right)\sum\limits_{n = 1}^\infty  {\left\{ {\sum\limits_{j = 2}^{p + 2m} {\frac{{{\zeta _n}\left( j \right)}}{{{n^{2p + 2m + 1 - j}}}}}  + 2\frac{{{H_n}}}{{{n^{2p + 2m}}}}} \right\}\left( {{x^n} - 1} \right)}  \\
&\quad  - {\left( { - 1} \right)^{p - 1}}\sum\limits_{n = 1}^\infty  {\left\{ {\sum\limits_{j = 2}^{p + 2m} {\frac{{{\zeta _n}\left( j \right)}}{{{n^{2p + 2m + 1 - j}}}}}  + 2\frac{{{H_n}}}{{{n^{2p + 2m}}}}} \right\}\left( {\sum\limits_{k = 1}^n {\frac{{{x^k}}}{k}} } \right)}   .\tag{2.14}
\end{align*}
After simplification, we find that
\begin{align*}
&{\left( { - 1} \right)^{p - 1}}\sum\limits_{n = 1}^\infty  {\left\{ {\sum\limits_{j = p + 1}^{p + 2m} {\frac{{{\zeta _n}\left( j \right)}}{{{n^{2p + 2m + 1 - j}}}}} } \right\}\left( {\sum\limits_{k = 1}^n {\frac{{{x^k}}}{k}} } \right)} \\
& =\left( {p + 2m + 1} \right){I_{p,p + 2m + 1}}\left( x \right) - \left( {p + 1} \right){I_{p + 1,p + 2m}}\left( x \right) \\
&\quad  +{\left( { - 1} \right)^{p - 1}}\ln \left( {1 - x} \right)\sum\limits_{n = 1}^\infty  {\left\{ {\sum\limits_{j = p + 1}^{p + 2m} {\frac{{{\zeta _n}\left( j \right)}}{{{n^{2p + 2m + 1 - j}}}}} } \right\}\left( {{x^n} - 1} \right)} \\
&\quad  + \sum\limits_{i = 1}^{p + 2m - 1} {{{\left( { - 1} \right)}^{i - 1}}L{i_{p + 2m + 1 - i}}\left( x \right)\sum\limits_{n = 1}^\infty  {\left\{ {\sum\limits_{j = 2}^p {\frac{{{\zeta _n}\left( j \right)}}{{{n^{p + 1 + i - j}}}}}  + 2\frac{{{H_n}}}{{{n^{p + i}}}}} \right\}{x^n}} }   \\
&\quad   - \sum\limits_{i = 1}^{p - 1} {{{\left( { - 1} \right)}^{i - 1}}L{i_{p + 1 - i}}\left( x \right)\sum\limits_{n = 1}^\infty  {\left\{ {\sum\limits_{j = 2}^{p + 2m} {\frac{{{\zeta _n}\left( j \right)}}{{{n^{p + 2m + 1 + i - j}}}}}  + 2\frac{{{H_n}}}{{{n^{p + 2m + i}}}}} \right\}{x^n}} }  .\tag{2.15}
\end{align*}
Replacing $p$ by $p-1$, $m$ by $m+1$ in (2.15), then combining with (2.15), we can obtain (2.12).\\
Let $x=-1$ and $x\rightarrow1$ in (2.12), we can gives the following Corollaries:
\begin{cor}
For $p\geq 2,m\geq 0$ and $p,m\in Z$, we have
\begin{align*}
&{\left( { - 1} \right)^p}\sum\limits_{n = 1}^\infty  {\left\{ {\frac{{{H_n}{\zeta _n}\left( {p + 2m + 1} \right)}}{{{n^p}}} + \frac{{{H_n}{\zeta _n}\left( p \right)}}{{{n^{p + 2m + 1}}}}} \right\}} \\
& =\left( {p + 2m + 2} \right)I_{p - 1,p + 2m + 2}(1) + \left( {2m + 1} \right)I_{p,p + 2m + 1}(1)- \left( {p + 1} \right)I_{p + 1,p + 2m}(1)\\
&\quad+ \sum\limits_{i = 1}^{p + 2m} {{{\left( { - 1} \right)}^{i - 1}}\zeta \left( {p + 2m + 2 - i} \right)\sum\limits_{n = 1}^\infty  {\left\{ {\sum\limits_{j = 2}^{p - 1} {\frac{{{\zeta _n}\left( j \right)}}{{{n^{p + i - j}}}}}  + 2\frac{{{H_n}}}{{{n^{p - 1 + i}}}}} \right\}} } \\
&\quad + \sum\limits_{i = 1}^{p + 2m - 1} {{{\left( { - 1} \right)}^{i - 1}}\zeta \left( {p + 2m + 1 - i} \right)\sum\limits_{n = 1}^\infty  {\left\{ {\sum\limits_{j = 2}^p {\frac{{{\zeta _n}\left( j \right)}}{{{n^{p + i + 1 - j}}}}}  + 2\frac{{{H_n}}}{{{n^{p + i}}}}} \right\}} } \\
&\quad - \sum\limits_{i = 1}^{p - 2} {{{\left( { - 1} \right)}^{i - 1}}\zeta \left( {p - i} \right)\sum\limits_{n = 1}^\infty  {\left\{ {\sum\limits_{j = 2}^{p + 2m + 1} {\frac{{{\zeta _n}\left( j \right)}}{{{n^{p + 2m + i + 2 - j}}}}}  + 2\frac{{{H_n}}}{{{n^{p + 2m + 1 + i}}}}} \right\}} }  \\
&\quad - \sum\limits_{i = 1}^{p - 1} {{{\left( { - 1} \right)}^{i - 1}}\zeta \left( {p + 1 - i} \right)\sum\limits_{n = 1}^\infty  {\left\{ {\sum\limits_{j = 2}^{p + 2m} {\frac{{{\zeta _n}\left( j \right)}}{{{n^{p + 2m + i + 1 - j}}}}}  + 2\frac{{{H_n}}}{{{n^{p + 2m + i}}}}} \right\}} } .\tag{2.16}
\end{align*}
\end{cor}
\begin{cor}
For $p\geq 2,m\geq 0$ and $p,m\in Z$, we have
\begin{align*}
&{\left( { - 1} \right)^p}\sum\limits_{n = 1}^\infty  {\left\{ {\frac{{{L_n}\left( 1 \right){\zeta _n}\left( {p + 2m + 1} \right)}}{{{n^p}}} + \frac{{{L_n}\left( 1 \right){\zeta _n}\left( p \right)}}{{{n^{p + 2m + 1}}}}} \right\}}  \\
& =\left( {p + 1} \right){I_{p + 1,p + 2m}}\left( { - 1} \right) - \left( {2m + 1} \right){I_{p,p + 2m + 1}}\left( { - 1} \right) - \left( {p + 2m + 2} \right){I_{p - 1,p + 2m + 2}}\left( { - 1} \right)\\
&\quad + {\left( { - 1} \right)^p}\ln 2\sum\limits_{n = 1}^\infty  {\left\{ {\frac{{{\zeta _n}\left( {p + 2m + 1} \right)}}{{{n^p}}} + \frac{{{\zeta _n}\left( p \right)}}{{{n^{p + 2m + 1}}}}} \right\}\left( {{{\left( { - 1} \right)}^{n - 1}} + 1} \right)}  \\
&\quad  + \sum\limits_{i = 1}^{p - 1} {{{\left( { - 1} \right)}^{i - 1}}\bar \zeta \left( {p + 1 - i} \right)\sum\limits_{n = 1}^\infty  {\left\{ {\sum\limits_{j = 2}^{p + 2m} {\frac{{{\zeta _n}\left( j \right)}}{{{n^{p + 2m + 1 + i - j}}}}}  + 2\frac{{{H_n}}}{{{n^{p + 2m + i}}}}} \right\}{{\left( { - 1} \right)}^{n - 1}}} }    \\
&\quad  + \sum\limits_{i = 1}^{p - 2} {{{\left( { - 1} \right)}^{i - 1}}\bar \zeta \left( {p - i} \right)\sum\limits_{n = 1}^\infty  {\left\{ {\sum\limits_{j = 2}^{p + 2m + 1} {\frac{{{\zeta _n}\left( j \right)}}{{{n^{p + 2m + 2 + i - j}}}}}  + 2\frac{{{H_n}}}{{{n^{p + 2m + 1 + i}}}}} \right\}{{\left( { - 1} \right)}^{n - 1}}} } \\
&\quad - \sum\limits_{i = 1}^{p + 2m - 1} {{{\left( { - 1} \right)}^{i - 1}}\bar \zeta \left( {p + 2m + 1 - i} \right)\sum\limits_{n = 1}^\infty  {\left\{ {\sum\limits_{j = 2}^p {\frac{{{\zeta _n}\left( j \right)}}{{{n^{p + 1 + i - j}}}}}  + 2\frac{{{H_n}}}{{{n^{p + i}}}}} \right\}{{\left( { - 1} \right)}^{n - 1}}} }\\
&\quad  - \sum\limits_{i = 1}^{p + 2m} {{{\left( { - 1} \right)}^{i - 1}}\bar \zeta \left( {p + 2m + 2 - i} \right)\sum\limits_{n = 1}^\infty  {\left\{ {\sum\limits_{j = 2}^{p - 1} {\frac{{{\zeta _n}\left( j \right)}}{{{n^{p + i - j}}}}}  + 2\frac{{{H_n}}}{{{n^{p + i - 1}}}}} \right\}{{\left( { - 1} \right)}^{n - 1}}} }.  \tag{2.17}
\end{align*}
\end{cor}
In the same manner, we obtain the following Theorem:
\begin{thm}
For $p\geq 2,m\geq 0$ and $p,m\in Z$, we have
\begin{align*}
&{\left( { - 1} \right)^p}\sum\limits_{n = 1}^\infty  {\left\{ {\frac{{{H_n}{\zeta _n}\left( {p + 2m + 2} \right)}}{{{n^p}}} - \frac{{{H_n}{\zeta _n}\left( p \right)}}{{{n^{p + 2m + 2}}}}} \right\}} \\
& =\sum\limits_{i = 1}^{p + 2m + 1} {{{\left( { - 1} \right)}^{i - 1}}\zeta \left( {p + 2m + 3 - i} \right)\sum\limits_{n = 1}^\infty  {\left\{ {\sum\limits_{j = 2}^{p - 1} {\frac{{{\zeta _n}\left( j \right)}}{{{n^{p + i - j}}}}}  + 2\frac{{{H_n}}}{{{n^{p - 1 + i}}}}} \right\}} } \\
&\quad  + \sum\limits_{i = 1}^{p + 2m} {{{\left( { - 1} \right)}^{i - 1}}\zeta \left( {p + 2m + 2 - i} \right)\sum\limits_{n = 1}^\infty  {\left\{ {\sum\limits_{j = 2}^p {\frac{{{\zeta _n}\left( j \right)}}{{{n^{p + i + 1 - j}}}}}  + 2\frac{{{H_n}}}{{{n^{p + i}}}}} \right\}} }  \\
&\quad - \sum\limits_{i = 1}^{p - 2} {{{\left( { - 1} \right)}^{i - 1}}\zeta \left( {p - i} \right)\sum\limits_{n = 1}^\infty  {\left\{ {\sum\limits_{j = 2}^{p + 2m + 2} {\frac{{{\zeta _n}\left( j \right)}}{{{n^{p + 2m + i + 3 - j}}}}}  + 2\frac{{{H_n}}}{{{n^{p + 2m + 2 + i}}}}} \right\}} }  \\
&\quad - \sum\limits_{i = 1}^{p - 1} {{{\left( { - 1} \right)}^{i - 1}}\zeta \left( {p + 1 - i} \right)\sum\limits_{n = 1}^\infty  {\left\{ {\sum\limits_{j = 2}^{p + 2m + 1} {\frac{{{\zeta _n}\left( j \right)}}{{{n^{p + 2m + i + 2 - j}}}}}  + 2\frac{{{H_n}}}{{{n^{p + 2m + i + 1}}}}} \right\}} } \\
&\quad+ \left( {p + 2m + 3} \right)I\left( {p - 1,p + 2m + 3} \right) + \left( {2m + 2} \right)I\left( {p,p + 2m + 2} \right)\\
&\quad - \left( {p + 1} \right)I\left( {p + 1,p + 2m + 1} \right) .\tag{2.18}
\end{align*}
\end{thm}
\pf  Similarly to the proof of Theorem 2.5, considering integral
$\int\limits_0^1 {\frac{{\ln \left( {1 - x} \right)L{i_p}\left( x \right)L{i_{p + 2m + 1}}\left( x \right)}}{x}} dx$ , we deduce Theorem 2.8 holds.
\begin{thm}
For $p\geq 2,m\geq 0$ and $p,m\in Z$, we have
\begin{align*}
&{\left( { - 1} \right)^p}\sum\limits_{n = 1}^\infty  {\left\{ {\frac{{{L_n}\left( 1 \right){L_n}\left( {p + 2m + 1} \right)}}{{{n^p}}} + \frac{{{L_n}\left( 1 \right){L_n}\left( p \right)}}{{{n^{p + 2m + 1}}}}} \right\}{{\left( { - 1} \right)}^n}}  \\
& =\left( {p + 2m + 1} \right){I_{p + 2m + 2,p - 1}}\left( { - 1} \right) + {R_{p + 2m + 2,p - 1}}\left( { - 1} \right) - \left( {p - 1} \right){I_{p,p + 2m + 1}}\left( { - 1} \right) - {R_{p,p + 2m + 1}}\left( { - 1} \right) \\
&\quad + \left( {p + 2m} \right){I_{p + 2m + 1,p}}\left( { - 1} \right) + {R_{p + 2m + 1,p}}\left( { - 1} \right) - p{I_{p + 1,p + 2m}}\left( { - 1} \right) - {R_{p + 1,p + 2m}}\left( { - 1} \right)
 \\
&\quad + {\left( { - 1} \right)^p}\ln 2\sum\limits_{n = 1}^\infty  {\left\{ {\frac{{{L_n}\left( {p + 2m + 1} \right)}}{{{n^p}}} + \frac{{{L_n}\left( p \right)}}{{{n^{p + 2m + 1}}}}} \right\}\left( {{{\left( { - 1} \right)}^n} - 1} \right)} \\
&\quad  + \sum\limits_{i = 1}^{p - 1} {{{\left( { - 1} \right)}^{i - 1}}\bar \zeta \left( {p + 1 - i} \right)\sum\limits_{n = 1}^\infty  {\left\{ {\frac{{{L_n}\left( 1 \right)}}{{{n^{p + 2m + i}}}}{{\left( { - 1} \right)}^{n - 1}} - \sum\limits_{j = 2}^{p + 2m} {\frac{{{L_n}\left( j \right)}}{{{n^{p + 2m + 1 + i - j}}}}} } \right\}} }   \\
&\quad + \sum\limits_{i = 1}^{p - 2} {{{\left( { - 1} \right)}^{i - 1}}\bar \zeta \left( {p - i} \right)\sum\limits_{n = 1}^\infty  {\left\{ {\frac{{{L_n}\left( 1 \right)}}{{{n^{p + 2m + 1 + i}}}}{{\left( { - 1} \right)}^{n - 1}} - \sum\limits_{j = 2}^{p + 2m + 1} {\frac{{{L_n}\left( j \right)}}{{{n^{p + 2m + 2 + i - j}}}}} } \right\}} }  \\
&\quad - \sum\limits_{i = 1}^{p + 2m} {{{\left( { - 1} \right)}^{i - 1}}\bar \zeta \left( {p + 2m + 2 - i} \right)\sum\limits_{n = 1}^\infty  {\left\{ {\frac{{{L_n}\left( 1 \right)}}{{{n^{p + i - 1}}}}{{\left( { - 1} \right)}^{n - 1}} - \sum\limits_{j = 2}^{p - 1} {\frac{{{L_n}\left( j \right)}}{{{n^{p + i - j}}}}} } \right\}} } \\
&\quad  - \sum\limits_{i = 1}^{p + 2m - 1} {{{\left( { - 1} \right)}^{i - 1}}\bar \zeta \left( {p + 2m + 1 - i} \right)\sum\limits_{n = 1}^\infty  {\left\{ {\frac{{{L_n}\left( 1 \right)}}{{{n^{p + i}}}}{{\left( { - 1} \right)}^{n - 1}} - \sum\limits_{j = 2}^p {\frac{{{L_n}\left( j \right)}}{{{n^{p + 1 + i - j}}}}} } \right\}} }  .\tag{2.19}
\end{align*}
where ${R_{p,q}}\left( x \right) = \int\limits_0^x {\frac{{L{i_p}\left( { - t} \right)L{i_q}\left( t \right)}}{t}dt} $ and
\begin{align*}
{R_{p,q}}\left( x \right) &= \sum\limits_{i = 1}^{q - 1} {{{\left( { - 1} \right)}^{i - 1}}L{i_{p + i}}\left( { - x} \right)L{i_{q + 1 - i}}\left( x \right)}  + {\left( { - 1} \right)^q}\ln \left( {1 - x} \right)\left( {L{i_{p + q}}\left( { - x} \right) + \bar \zeta \left( {p + q} \right)} \right)
\nonumber \\ &\quad  - {\left( { - 1} \right)^q}\sum\limits_{n = 1}^\infty  {\frac{{{{\left( { - 1} \right)}^n}}}{{{n^{p + q}}}}\left( {\sum\limits_{k = 1}^n {\frac{{{x^k}}}{k}} } \right)}  . \tag{2.20}
\end{align*}
For example
\begin{align*}
&{R_{4,1}}\left( { - 1} \right) = \sum\limits_{n = 1}^\infty  {\frac{{{L_n}\left( 1 \right)}}{{{n^5}}}{{\left( { - 1} \right)}^{n - 1}}}  - \frac{{31}}{{16}}\zeta \left( 5 \right)\ln 2,\\
&{R_{2,3}}\left( { - 1} \right) = \sum\limits_{n = 1}^\infty  {\frac{{{L_n}\left( 1 \right)}}{{{n^5}}}{{\left( { - 1} \right)}^{n - 1}}}  + \frac{7}{8}\zeta \left( 6 \right) - \frac{3}{4}{\zeta ^2}\left( 3 \right) - \frac{{31}}{{16}}\zeta \left( 5 \right)\ln 2.
\end{align*}
\end{thm}
\pf  Similarly to the proof of Theorem 2.5 and 2.8, we consider the following integral
$$\int\limits_0^{ - 1} {\frac{{\ln \left( {1 + x} \right)L{i_p}\left( x \right)L{i_{p + 2m}}\left( x \right)}}{x}} dx.$$
Noting that
\[\ln \left( {1 + x} \right)L{i_p}\left( x \right) = pL{i_{p + 1}}\left( x \right) + L{i_{p + 1}}\left( { - x} \right) + \sum\limits_{n = 1}^\infty  {\frac{{{L_n}\left( 1 \right)}}{{{n^p}}}{x^n}}  + \sum\limits_{i = 1}^p {\sum\limits_{n = 1}^\infty  {\frac{{{L_n}\left( i \right)}}{{{n^{p + 1 - i}}}}{{\left( { - x} \right)}^n}} } . \tag{2.21}\]
After simplification, we deduce Theorem 2.9 holds.
\section{Representation of Euler sums by zeta values and linear sums}
In this section, we consider the analytic representations of two quadratic Euler sums which involves
harmonic numbers and alternating harmonic numbers
\[\sum\limits_{n = 1}^\infty  {\frac{{{\zeta _n}\left( p \right){\zeta _n}\left( {p + 2m + 1} \right)}}{n}} {\left( { - 1} \right)^{n - 1}},\sum\limits_{n = 1}^\infty  {\frac{{{L_n}\left( p \right){L_n}\left( {p + 2m + 1} \right)}}{n}} {\left( { - 1} \right)^{n - 1}}\]through zeta values and linear sums, and give explicit formulae for several 6-order quadratic sums in terms of zeta values and linear sums.
\begin{thm}
For $1 \le {l_1},{l_2},m \in Z$ and $x,y,x \in \left[ { - 1,1} \right)$, we have the following relation
\begin{align*}
&\sum\limits_{n = 1}^\infty  {\frac{{{\zeta _n}\left( {{l_1},x} \right){\zeta _n}\left( {{l_2},y} \right)}}{{{n^m}}}{z^n}}  + \sum\limits_{n = 1}^\infty  {\frac{{{\zeta _n}\left( {{l_1},x} \right){\zeta _n}\left( {m,z} \right)}}{{{n^{{l_2}}}}}{y^n}}  + \sum\limits_{n = 1}^\infty  {\frac{{{\zeta _n}\left( {{l_2},y} \right){\zeta _n}\left( {m,z} \right)}}{{{n^{{l_1}}}}}{x^n}} \\
& =\sum\limits_{n = 1}^\infty  {\frac{{{\zeta _n}\left( {m,z} \right)}}{{{n^{{l_1} + {l_2}}}}}{{\left( {xy} \right)}^n}}  + \sum\limits_{n = 1}^\infty  {\frac{{{\zeta _n}\left( {{l_1},x} \right)}}{{{n^{m + {l_2}}}}}{{\left( {yz} \right)}^n}}  + \sum\limits_{n = 1}^\infty  {\frac{{{\zeta _n}\left( {{l_2},y} \right)}}{{{n^{{l_1} + m}}}}{{\left( {xz} \right)}^n}}  \\
&\quad  +L{i_m}\left( z \right)L{i_{{l_1}}}\left( x \right)L{i_{{l_2}}}\left( y \right) - L{i_{{l_1} + {l_2} + m}}\left( {xyz} \right) \tag{3.1}
\end{align*}
where the partial sum ${\zeta _n}\left( {l,x} \right)$ is defined by ${\zeta _n}\left( {l,x} \right) = \sum\limits_{k = 1}^n {\frac{{{x^k}}}{k^l}}$.
\end{thm}
\pf First, we construct the function $F\left( {x,y,z} \right) = \sum\limits_{n = 1}^\infty  {\left\{ {{\zeta _n}\left( {{l_1},x} \right){\zeta _n}\left( {{l_2},y} \right) - {\zeta _n}\left( {{l_1} + {l_2},xy} \right)} \right\}{z^{n - 1}}} $. By the definition of ${\zeta _n}\left( {l,x} \right)$, we have
\[F\left( {x,y,z} \right) = zF\left( {x,y,z} \right) + \sum\limits_{n = 1}^\infty  {\left\{ {\frac{{{\zeta _n}\left( {{l_1},x} \right)}}{{{{\left( {n + 1} \right)}^{{l_2}}}}}{y^{n+1}} + \frac{{{\zeta _n}\left( {{l_2},y} \right)}}{{{{\left( {n + 1} \right)}^{{l_1}}}}}{x^{n+1}}} \right\}{z^n}} .\tag{3.2}\]
Moving $zF(x,y,z)$ from right to left and then multiplying  $(1-z)^{-1}$  to the equation (3.2) and integrating over the interval $(0,z)$, we obtain
\[\sum\limits_{n = 1}^\infty  {\frac{{{\zeta _n}\left( {{l_1},x} \right){\zeta _n}\left( {{l_2},y} \right) - {\zeta _n}\left( {{l_1} + {l_2},xy} \right)}}{n}{z^{n}}}  = \sum\limits_{n = 1}^\infty  {\left\{ {\frac{{{\zeta _n}\left( {{l_1},x} \right)}}{{{{\left( {n + 1} \right)}^{{l_2}}}}}{y^{n+1}} + \frac{{{\zeta _n}\left( {{l_2},y} \right)}}{{{{\left( {n + 1} \right)}^{{l_1}}}}}{x^{n+1}}} \right\}\left\{ {L{i_1}\left( z \right) - {\zeta _n}\left( {1,z} \right)} \right\}} \tag{3.3}\]
Furthermore, using integration and the following formula
\[\sum\limits_{n = 1}^\infty  {\left\{ {\frac{{{\zeta _n}\left( {{l_1},x} \right)}}{{{{\left( {n + 1} \right)}^{{l_2}}}}}{y^{n + 1}} + \frac{{{\zeta _n}\left( {{l_2},y} \right)}}{{{{\left( {n + 1} \right)}^{{l_1}}}}}{x^{n + 1}}} \right\}}  = L{i_{{l_1}}}\left( x \right)L{i_{{l_2}}}\left( y \right) - L{i_{{l_1} + {l_2}}}\left( {xy} \right)\]
we can obtain (3.1). \\Let $(x,y,z)=(-1,1,1), (l_1,l_2,m)=(1,p+2m+1,p)$ and $(x,y,z)=(-1,-1,-1), (l_1,l_2,m)=(1,p+2m+1,p)$ in (3.1), we can gives the following Corollaries
\begin{cor}
If $1<p\in Z, 0\leq m\in Z$, then we have
\begin{align*}
&\sum\limits_{n = 1}^\infty  {\frac{{{L_n}\left( 1 \right){\zeta _n}\left( {p + 2m + 1} \right)}}{{{n^p}}}}  + \sum\limits_{n = 1}^\infty  {\frac{{{L_n}\left( 1 \right){\zeta _n}\left( p \right)}}{{{n^{p + 2m + 1}}}}}  + \sum\limits_{n = 1}^\infty  {\frac{{{\zeta _n}\left( {p + 2m + 1} \right){\zeta _n}\left( p \right)}}{n}{{\left( { - 1} \right)}^{n - 1}}}  \\
& =\sum\limits_{n = 1}^\infty  {\frac{{{\zeta _n}\left( p \right)}}{{{n^{p + 2m + 2}}}}{{\left( { - 1} \right)}^{n - 1}}}  + \sum\limits_{n = 1}^\infty  {\frac{{{L_n}\left( 1 \right)}}{{{n^{2p + 2m + 1}}}}}  + \sum\limits_{n = 1}^\infty  {\frac{{{\zeta _n}\left( {p + 2m + 1} \right)}}{{{n^{p + 1}}}}{{\left( { - 1} \right)}^{n - 1}}}  \\
&\quad  +\ln 2\zeta \left( {p + 2m + 1} \right)\zeta \left( p \right) - \bar \zeta \left( {2p + 2m + 2} \right). \tag{3.4}
\end{align*}
\end{cor}
\begin{cor}
If $1\leq p\in Z, 0\leq m\in Z$, then we have
\begin{align*}
&\sum\limits_{n = 1}^\infty  {\frac{{{L_n}\left( 1 \right){L_n}\left( {p + 2m + 1} \right)}}{{{n^p}}}{{\left( { - 1} \right)}^{n - 1}}}  + \sum\limits_{n = 1}^\infty  {\frac{{{L_n}\left( 1 \right){L_n}\left( p \right)}}{{{n^{p + 2m + 1}}}}{{\left( { - 1} \right)}^{n - 1}}}  + \sum\limits_{n = 1}^\infty  {\frac{{{L_n}\left( {p + 2m + 1} \right){L_n}\left( p \right)}}{n}{{\left( { - 1} \right)}^{n - 1}}}  \\
& =\sum\limits_{n = 1}^\infty  {\frac{{{L_n}\left( p \right)}}{{{n^{p + 2m + 2}}}}}  + \sum\limits_{n = 1}^\infty  {\frac{{{L_n}\left( 1 \right)}}{{{n^{2p + 2m + 1}}}}}  + \sum\limits_{n = 1}^\infty  {\frac{{{L_n}\left( {p + 2m + 1} \right)}}{{{n^{p + 1}}}}}  \\
&\quad  +\ln 2\bar \zeta \left( {p + 2m + 1} \right)\bar \zeta \left( p \right) - \bar \zeta \left( {2p + 2m + 2} \right). \tag{3.5}
\end{align*}
\end{cor}
From Corollary 2.7, 3.2, 3.3 and Theorem 2.9, we can obtain the following Theorem
\begin{thm} For $2\leq p \in Z$ and $0\leq m \in Z$, the quadratic sums
\[\sum\limits_{n = 1}^\infty  {\frac{{{\zeta _n}\left( p \right){\zeta _n}\left( {p + 2m + 1} \right)}}{n}} {\left( { - 1} \right)^{n - 1}},\sum\limits_{n = 1}^\infty  {\frac{{{L_n}\left( p \right){L_n}\left( {p + 2m + 1} \right)}}{n}} {\left( { - 1} \right)^{n - 1}}\]
are reducible to linear sums.
\end{thm}
Using mathematica, we can obtain the following numerical values
\begin{align*}
&{{L}}{{{i}}_4}\left( {\frac{1}{2}} \right)= 0.5174790616738993863307581618988629456,\\
&\sum\limits_{n = 1}^\infty  {\frac{{{L_n}\left( 1 \right)}}{{{n^5}}}{{\left( { - 1} \right)}^{n - 1}}}  = 0.987441426403299713771650007985,\\
&\sum\limits_{n = 1}^\infty  {\frac{{{L_n}\left( 2 \right)}}{{{n^4}}}}  = 1.06358224101814909880154833539,\\
&\sum\limits_{n = 1}^\infty  {\frac{{{\zeta _n}\left( 2 \right)}}{{{n^4}}}{{\left( { - 1} \right)}^{n - 1}}}  = 0.934707899349253255197542851216,\\
&\sum\limits_{n = 1}^\infty  {\frac{{{H_n}}}{{{n^5}}}{{\left( { - 1} \right)}^{n - 1}}}  = 0.959151942504318157165421137321,\\
&\sum\limits_{n = 1}^\infty  {\frac{{{L_n}\left( 1 \right)}}{{{n^5}}}}  = 1.02005194570145237930331996837,\\
&\sum\limits_{n = 1}^\infty  {\frac{{{L_n}\left( 1 \right){\zeta _n}\left( 2 \right)}}{{{n^3}}}}  = 1.15935334356951415975457027807,\\
&\sum\limits_{n = 1}^\infty  {\frac{{{L_n}\left( 1 \right){\zeta _n}\left( 3 \right)}}{{{n^2}}}}  = 1.47723102170162037670053143416,\\
\end{align*}
Taking $p=2,m=0$ in (2.17)(2.19)(3.4)(3.5), we obtain
\begin{align*}
&\sum\limits_{n = 1}^\infty  {\frac{{{L_n}\left( 1 \right){\zeta _n}\left( 3 \right)}}{{{n^2}}}}  + \sum\limits_{n = 1}^\infty  {\frac{{{L_n}\left( 1 \right){\zeta _n}\left( 2 \right)}}{{{n^3}}}}  \\
& =\frac{3}{4}{\zeta ^2}\left( 3 \right) + \frac{7}{4}\zeta \left( 6 \right) + \frac{5}{8}\zeta \left( 2 \right)\zeta \left( 3 \right)\ln 2 - 2\zeta \left( 2 \right)L{i_4}\left( {\frac{1}{2}} \right) + \frac{5}{4}\zeta \left( 4 \right){\ln ^2}2 \\ &\quad- \frac{1}{{12}}\zeta \left( 2 \right){\ln ^4}2, \tag{3.6}
\end{align*}
\begin{align*}
&\sum\limits_{n = 1}^\infty  {\frac{{{\zeta _n}\left( 2 \right){\zeta _n}\left( 3 \right)}}{n}{{\left( { - 1} \right)}^{n - 1}}}  \\
& = - \frac{{161}}{{64}}\zeta \left( 6 \right) + \frac{{31}}{{16}}\zeta \left( 5 \right)\ln 2 + \frac{9}{{32}}{\zeta ^2}\left( 3 \right) + \frac{3}{8}\zeta \left( 2 \right)\zeta \left( 3 \right)\ln 2 + 2\zeta \left( 2 \right)L{i_4}\left( {\frac{1}{2}} \right) \\
&\quad  - \frac{5}{4}\zeta \left( 4 \right){\ln ^2}2 + \frac{1}{{12}}\zeta \left( 2 \right){\ln ^4}2 + \sum\limits_{n = 1}^\infty  {\frac{{{\zeta _n}\left( 2 \right)}}{{{n^4}}}{{\left( { - 1} \right)}^{n - 1}}}  - \sum\limits_{n = 1}^\infty  {\frac{{{L_n}\left( 3 \right)}}{{{n^3}}}} , \tag{3.7}
\end{align*}
\begin{align*}
&\sum\limits_{n = 1}^\infty  {\frac{{{L_n}\left( 1 \right){L_n}\left( 3 \right)}}{{{n^2}}}{{\left( { - 1} \right)}^{n - 1}}}  + \sum\limits_{n = 1}^\infty  {\frac{{{L_n}\left( 1 \right){L_n}\left( 2 \right)}}{{{n^3}}}{{\left( { - 1} \right)}^{n - 1}}}   \\
& = - \frac{{385}}{{128}}\zeta \left( 6 \right) + \frac{{31}}{8}\zeta \left( 5 \right)\ln 2 + \frac{3}{{32}}{\zeta ^2}\left( 3 \right) + \frac{9}{8}\zeta \left( 2 \right)\zeta \left( 3 \right)\ln 2 + \zeta \left( 2 \right)L{i_4}\left( {\frac{1}{2}} \right) \\
&\quad - \frac{5}{8}\zeta \left( 4 \right){\ln ^2}2 + \frac{1}{{24}}\zeta \left( 2 \right){\ln ^4}2\, \tag{3.8}
\end{align*}
\begin{align*}
&\sum\limits_{n = 1}^\infty  {\frac{{{L_n}\left( 2 \right){L_n}\left( 3 \right)}}{n}{{\left( { - 1} \right)}^{n - 1}}}  \\
& = \frac{{163}}{{128}}\zeta \left( 6 \right) - \frac{{31}}{{16}}\zeta \left( 5 \right)\ln 2 + \frac{3}{{16}}{\zeta ^2}\left( 3 \right) - \frac{3}{4}\zeta \left( 2 \right)\zeta \left( 3 \right)\ln 2 - \zeta \left( 2 \right)L{i_4}\left( {\frac{1}{2}} \right)\\
&\quad + \frac{5}{8}\zeta \left( 4 \right){\ln ^2}2 - \frac{1}{{24}}\zeta \left( 2 \right){\ln ^4}2 + \sum\limits_{n = 1}^\infty  {\frac{{{L_n}\left( 2 \right)}}{{{n^4}}}}  + \sum\limits_{n = 1}^\infty  {\frac{{{L_n}\left( 3 \right)}}{{{n^3}}}} . \tag{3.9}
\end{align*}
In [21], we gave the following formula
\begin{align*}
\sum\limits_{n = 1}^\infty  {\frac{{{L_n}\left( 1 \right){\zeta _n}\left( 2 \right)}}{{{n^3}}}} &=\frac{{29}}{8}\zeta \left( 2 \right)\zeta \left( 3 \right)\ln 2 - \frac{{93}}{{32}}\zeta \left( 5 \right)\ln 2 - \frac{{1855}}{{128}}\zeta \left( 6 \right) + \frac{{17}}{{16}}{\zeta ^2}\left( 3 \right)- \sum\limits_{n = 1}^\infty  {\frac{{{L_n}\left( 1 \right)}}{{{n^5}}}{{\left( { - 1} \right)}^{n - 1}}} \\
& \quad + \sum\limits_{n = 1}^\infty  {\frac{{{L_n}\left( 2 \right)}}{{{n^4}}}}  + 4\sum\limits_{n = 1}^\infty  {\frac{{{\zeta _n}\left( 2 \right)}}{{{n^4}}}{{\left( { - 1} \right)}^{n - 1}}}  + 8\sum\limits_{n = 1}^\infty  {\frac{{{H_n}}}{{{n^5}}}{{\left( { - 1} \right)}^{n - 1}}}.\tag{3.10}
\end{align*}
Substituting (3.10) into (3.6) respectively, we obtain
\begin{align*}
\sum\limits_{n = 1}^\infty  {\frac{{{L_n}\left( 1 \right){\zeta _n}\left( 3 \right)}}{{{n^2}}}}  &=\frac{{2079}}{{128}}\zeta \left( 6 \right) + \frac{{93}}{{32}}\zeta \left( 5 \right)\ln 2 - \frac{5}{{16}}{\zeta ^2}\left( 3 \right) - 3\zeta \left( 2 \right)\zeta \left( 3 \right)\ln 2 - 2\zeta \left( 2 \right)L{i_4}\left( {\frac{1}{2}} \right) \\
& \quad + \frac{5}{4}\zeta \left( 4 \right){\ln ^2}2 - \frac{1}{{12}}\zeta \left( 2 \right){\ln ^4}2 + \sum\limits_{n = 1}^\infty  {\frac{{{L_n}\left( 1 \right)}}{{{n^5}}}{{\left( { - 1} \right)}^{n - 1}}}  - \sum\limits_{n = 1}^\infty  {\frac{{{L_n}\left( 2 \right)}}{{{n^4}}}}\\
& \quad - 4\sum\limits_{n = 1}^\infty  {\frac{{{\zeta _n}\left( 2 \right)}}{{{n^4}}}{{\left( { - 1} \right)}^{n - 1}}}  - 8\sum\limits_{n = 1}^\infty  {\frac{{{H_n}}}{{{n^5}}}{{\left( { - 1} \right)}^{n - 1}}} .\tag{3.11}
\end{align*}
Proceeding in a similar fashion to evaluation of the Theorem 2.5, it is possible to evaluate other Euler sums involving harmonic numbers and alternating harmonic numbers. For example
\[\sum\limits_{n = 1}^\infty  {\left\{ {\frac{{{H_n}{L_n}\left( 1 \right)}}{{{n^3}}} - \frac{{{H_n}{L_n}\left( 3 \right)}}{n}} \right\}{{\left( { - 1} \right)}^{n - 1}}}  = \frac{{15}}{4}\zeta \left( 4 \right)\ln 2 - \frac{9}{8}\zeta \left( 2 \right)\zeta \left( 3 \right) - \frac{1}{2}\zeta \left( 3 \right){\ln ^2}2.\]
{\bf Acknowledgments.} The authors would like to thank the anonymous
referee for his/her helpful comments, which improve the presentation
of the paper.

 \end{document}